\documentclass[a4paper]{amsart}

\usepackage[dvips]{graphicx}

\newtheorem{theorem}{Theorem}[section]

\newtheorem{lemma}[theorem]{Lemma}
\newtheorem{claim}[theorem]{Claim}

\newtheorem{corollary}[theorem]{Corollary}

\theoremstyle{definition}

\newtheorem{example}[theorem]{Example}

\theoremstyle{remark}

\numberwithin{equation}{section}

\begin{document}

\title[Crosscap numbers of 2-bridge knots]
{Crosscap numbers of 2-bridge knots}
\author[M. Hirasawa]{Mikami Hirasawa}
\address{Department of Mathematics, Gakushuin University,
1-5-1 Mejiro, Toshima-ku, Tokyo 171-8588, Japan}
\thanks{The first author was partially supported by MEXT, Grant-in-Aid for Young Scientists (B) 14740048.}
\email{hirasawa@math.gakushuin.ac.jp}
\author[M. Teragaito]{Masakazu Teragaito}
\address{Department of Mathematics and Mathematics Education,
Faculty of Education, 
Hiroshima University, 1-1-1 Kagamiyama, Higashi-hiroshima 739-8524, Japan}
\email{teragai@hiroshima-u.ac.jp}



\begin{abstract}
We present a practical algorithm to determine
the minimal genus of non-orientable spanning surfaces for $2$-bridge knots, called
the crosscap numbers.
We will exhibit a table of crosscap numbers of $2$-bridge knots up to $12$ crossings (all $362$ of them).
\end{abstract}
\maketitle

\section{Introduction}\label{sec:intro}

For a knot $K$ in the $3$-sphere $S^3$, there is a connected compact
embedded surface $F$ in $S^3$
whose boundary is $K$.
In particular, $F$ can be chosen to be orientable, and then it is called a Seifert surface for $K$.
The genus $g(K)$ of $K$ is the minimal number of genera of all Seifert surfaces for $K$.
Thus the unknot is the only knot of genus zero.

On the other hand, we can choose the above $F$ to be non-orientable,
for example, by adding a half-twisted band to a Seifert surface.
In this paper, such $F$ is referred to as a \textit{non-orientable spanning surface\/} for $K$.
We define the \textit{crosscap number\/} $\gamma(K)$ of a non-trivial knot $K$
as the minimal number of the first betti numbers $\beta_1$ of all non-orientable spanning surfaces for $K$, and
set $\gamma(\mathrm{unknot})=0$ for convenience.
We call $\gamma(K)$ the crosscap number because it counts the number of
`crosscap summands' in the closed surface obtained by capping off a non-orientable spanning surface
with a disk, which is well known to be a connected sum of projective planes.
In the literature, a crosscap number is also called a non-orientable genus \cite{MY2}.
For a non-trivial knot $K$, if a non-orientable spanning surface $F$ satisfies $\beta_1(F)=\gamma(K)$,
then $F$ is called a \textit{minimal genus non-orientable spanning surface\/} for $K$.

In general, it is very hard to determine the crosscap number for a given knot.
Any minimal genus Seifert surface becomes a non-orientable spanning surface for the same knot
if we attach a small half-twisted band as above, and hence we have an obvious inequality
$\gamma(K)\le 2g(K)+1$.
There are only a few results about crosscap numbers of knots.
Clark \cite{C} introduced the notion of crosscap number and
pointed out that $\gamma(K)=1$ if and only if $K$ is a $2$-cabled knot.
He also asked the existence of a knot satisfying the equality $\gamma(K)=2g(K)+1$, and
Murakami and Yasuhara \cite{MY} came up with the first example,
showing $\gamma(7_4)=3$ algebraically.
In \cite{T1}, the crosscap numbers of torus knots are completely determined.

The purpose of this paper is to determine the crosscap numbers of $2$-bridge knots, which form
a special but important class of knots.
For $2$-bridge knots, Hatcher and Thurston \cite{HT} 
constructed all incompressible, boundary-incompressible orientable or non-orientable spanning surfaces.
However, for the $2$-bridge knot $7_4$, a minimal genus non-orientable spanning surface can be
realized only by a boundary-compressible surface.
Then Bessho \cite{B} proved that any incompressible, boundary-compressible spanning surface
for a $2$-bridge knot becomes an incompressible, boundary-incompressible surface after
several boundary-compressions.
Therefore, theoretically, we can obtain $\gamma(K)$ as follows:
\begin{quote}
For a $2$-bridge knot $K$, generate all incompressible, boundary-incompressible spanning
surfaces according to \cite{HT}.
Let $n$ be the minimal first betti number of them.
Then if $n$ is realized by a non-orientable spanning surface, then $\gamma(K)=n$, and otherwise
$\gamma(K)=n+1$.
Here, $n$ equals the minimal length of all continued fraction expansions for $K$.
\end{quote}
However, an effective algorithm to determine $n$ was missing, and one could not tell,
for example, for which $2$-bridge knots, the equality $\gamma(K)=2g(K)+1$ holds.

In the following, we present a practical algorithm to find a shortest
continued fraction expansion for all rational numbers representing a $2$-bridge knot $K$.
This enables us to determine the crosscap number from any continued fraction expansion for $K$.
The main tool is so-called the modular diagram, whose vertices correspond to rational numbers,
on which we introduce the notion of depth.
In Section \ref{sec:table}, we exhibit a table of crosscap numbers of $2$-bridge knots up to $12$ crossings (all $362$ of them).

\section{Statement of results}\label{sec:result}


Let $K$ be a $2$-bridge knot $S(q,p)$ in Schubert's notation.
Here, $p$ and $q$ are coprime integers, and $q$ is odd.
As is well-known, $S(q,p)$ and $S(q',p')$ are equivalent if and only if $q'=q$ and $p'\equiv p^{\pm 1} \pmod q$, and
$S(q,-p)$ gives the mirror image of $S(q,p)$.

Consider a \textit{subtractive\/} continued fraction expansion of $p/q$ (see \cite{HT})

$$\frac{p}{q}=r+[b_1,b_2,\dots,b_n]\\
      =r+ \cfrac{1}{b_1
          -\cfrac{1}{b_2
          -\cfrac{1}{b_3
          -\cfrac{1}{\ddots
          -\cfrac{1}{b_n}}}}}$$
where $r, b_i\in \mathbb{Z}$ and $b_i\ne 0$.
The \textit{length\/} of this expansion is $n$.
Then $K$ is the boundary of the surface obtained by plumbing $n$ bands in a row, the
$i$th band having $b_i$ half-twists (right-handed if $b_i>0$ and left-handed if $b_i<0$).
If some $b_i$ is odd, then the expansion is said to be
of \textit{odd type}.
Otherwise, it is of \textit{even type}.
Any fraction has expansions of odd type and even type, e.g., $1/3=1-2/3=1+[-2,-2]=1+[-1,2]$.
In this paper, an expansion always means a subtractive one.
We remark the following equality:
$$r+[a_1,-a_2,a_3,-a_4,\dots,(-1)^{n-1}a_n]\\
      =r+  \cfrac{1}{a_1
          +\cfrac{1}{a_2
          +\cfrac{1}{a_3
          +\cfrac{1}{\ddots
          +\cfrac{1}{a_n}}}}}.$$

The crosscap number of a $2$-bridge knot $K$ can be described
in terms of the length of expansion corresponding to $K$.
The first theorem is due to Bessho, but we will give its proof
for reader's convenience in Section \ref{sec:bessho}.

\begin{figure}[ht]
\includegraphics*[scale=0.6]{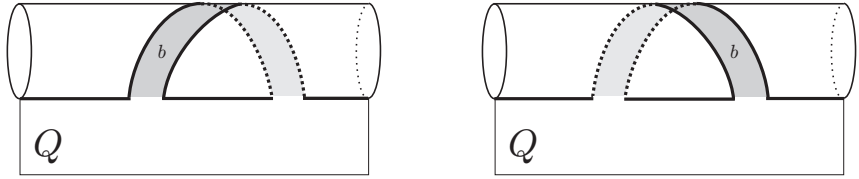}
\caption{M\"{o}bius bands attached to Seifert surfaces}\label{band}
\end{figure}

\begin{theorem}[Bessho \cite{B}]\label{thm:bessho}
Let $K$ be a $2$-bridge knot.
\begin{itemize}
\item[(1)] The crosscap number $\gamma(K)$ equals the minimal length
of  all expansions of odd type of all fractions corresponding to $K$.
\item[(2)] If a minimal genus non-orientable spanning surface $F$ for $K$ is boundary-compressible, then
$F$ is obtained from a minimal genus Seifert surface for $K$ by attaching a M\"{o}bius band as in Figure \ref{band}.
\end{itemize}
\end{theorem}

We present a practical algorithm to obtain a shortest expansion from any one of $p/q$.

\begin{theorem}\label{thm:algorithm}
Let $p/q=r+[b_1,b_2,\dots,b_n]$ be an expansion obtained from an arbitrary expansion of $p/q$
by fully reducing the length by a repetition of the following three reductions.
Then $n$ is the minimal length of all expansions of $p/q$.
\begin{itemize}
\item[(1)] Removal of coefficient $0$.
\[
\begin{array}{rcl}
[\dots,a,0,b,\dots] &=& [\dots,a+b,\dots],\\
\,[\dots,a,b,0] &=& [\dots,a],\\
\,[0,a,b,\dots] &=& -a+[b,\dots].\\
\end{array}
\]
\item[(2)] Removal of coefficient $\varepsilon=\pm 1$.
\[
\begin{array}{rcl}
[\dots,a,\varepsilon,b,\dots] &=& [\dots,a-\varepsilon,b-\varepsilon,\dots],\\
\,[\dots,a,\varepsilon] &=& [\dots,a-\varepsilon],\\
r+[\varepsilon,a,\dots] &=& (r+\varepsilon)+[a-\varepsilon,\dots].\\
\end{array}
\]
\item[(3)] Removal of a subsequence $2\varepsilon,2\varepsilon$ or $2\varepsilon,3\varepsilon,\dots,3\varepsilon,2\varepsilon$.
\textup{(}Here, $\varepsilon=\pm 1$, and possibly $m=2$.\textup{)}
\[
\begin{array}{rcl}
[\dots,a,\underbrace{2\varepsilon,3\varepsilon,\dots,3\varepsilon,2\varepsilon}_{m},b,\dots,] &=& [\dots,a-\varepsilon,\underbrace{-3\varepsilon,-3\varepsilon,\dots,-3\varepsilon}_{m-1},b-\varepsilon,\dots],\\
\,[\dots,a,\underbrace{2\varepsilon,3\varepsilon,\dots,3\varepsilon,2\varepsilon}_{m}] &=& [\dots,a-\varepsilon,\underbrace{-3\varepsilon,-3\varepsilon,\dots,-3\varepsilon}_{m-1}],\\
r+[\underbrace{2\varepsilon,3\varepsilon,\dots,3\varepsilon,2\varepsilon}_{m},a,\dots] &=& (r+\varepsilon)+[\underbrace{-3\varepsilon,-3\varepsilon,\dots,-3\varepsilon}_{m-1},a-\varepsilon,\dots],\\
r+[\underbrace{2\varepsilon,3\varepsilon,\dots,3\varepsilon,2\varepsilon}_{m}] &=& (r+\varepsilon)+[\underbrace{-3\varepsilon,-3\varepsilon,\dots,-3\varepsilon}_{m-1}].
\end{array}
\]
\end{itemize}
\end{theorem}

In Theorem \ref{thm:algorithm}, we fix a fraction $p/q$.
Although there are infinitely many fractions corresponding to a $2$-bridge knot, 
the next theorem guarantees that we can start from any fraction.

\begin{theorem}\label{thm:main}
Let $K$ be a $2$-bridge knot.
If the reduction in Theorem \ref{thm:algorithm} yields a length $n$ expansion,
then $n$ is the minimal length of all expansions of all fractions corresponding to $K$.
\end{theorem}

The next theorem is the key to determine whether a fraction $p/q$ admits a shortest expansion of odd type. 

\begin{theorem}\label{thm:shortest}
Two shortest expansions for $p/q$ are deformed to each other by
a finite repetition of the following, where $\varepsilon=\pm 1$\textup{:}
\[
\begin{array}{rcl}
[\dots,a,2\varepsilon,b,\dots] &=& [\dots,a-\varepsilon,-2\varepsilon,b-\varepsilon,\dots],\\
\,[\dots,a,2\varepsilon] &=& [\dots,a-\varepsilon,-2\varepsilon],\\
r+[2\varepsilon,a,\dots] &=& (r+\varepsilon)+[-2\varepsilon,a-\varepsilon,\dots].\\
\end{array}
\]
\end{theorem}


\begin{theorem}\label{thm:final}
Let $K=S(q,p)$ be a $2$-bridge knot.
If a shortest expansion of $p/q$ obtained by Theorem \ref{thm:algorithm}
contains an odd coefficient or $\pm 2$, then $\gamma(K)=n$, otherwise $\gamma(K)=n+1$, where
$n$ is the length of the expansion.
\end{theorem}

Remark that if there is a coefficient $\pm 2$ in an expansion with only even coefficients, then we can apply Theorem \ref{thm:shortest}
to obtain an expansion of odd type.


\begin{example}
Let $K=6_1$ in the knot table (see \cite{Ro}).
It is the $2$-bridge knot $S(9,2)$.
Then $2/9=[5,2]$. 
Thus $\gamma(K)=2$ by Theorem \ref{thm:final}.
\end{example}

It is known that any $2$-bridge knot $K$ has a unique expansion of even type modulo integer parts, and the length of which equals $2g(K)$.
As a direct corollary to Theorem \ref{thm:final}, we can completely characterize those $2$-bridge knots satisfying the equality
$\gamma(K)=2g(K)+1$.

\begin{corollary}\label{thm:2g+1}
For a $2$-bridge knot $K$, the equality $\gamma(K)=2g(K)+1$ holds if and only if
there is no coefficient $\pm 2$ in the \textup{(}unique\textup{)} expansion for $K$ containing only even coefficients.
\end{corollary}

\begin{example}\label{ex:2g+1}
Let $K=7_4=S(15,4)$.  Note that $K$ has genus one.
Since $4/15=[4,4]$, $\gamma(K)=2g(K)+1=3$ by Corollary \ref{thm:2g+1}.
More examples will be given in Section \ref{sec:example}.
\end{example}

Some minimal genus non-orientable surfaces for $2$-bridge knots are boundary-incompressible,
but others boundary-compressible, and some $2$-bridge knots have several such surfaces.
This makes a strong contrast to the case of torus knots, 
where minimal genus non-orientable spanning surfaces are 
boundary-incompressible and even unique \cite{T1}.

By the theorems above, we can characterize $2$-bridge knots with boundary-compressible minimal genus non-orientable spanning
surfaces.
It is unknown whether Corollary \ref{cor:compare} below generalizes to all knots.

\begin{theorem}\label{thm:characterization}
Let $K=S(q,p)$ be a $2$-bridge knot, and $\mathcal{C}$ the set of shortest expansions for $p/q$.
Then we have\textup{:}
\begin{itemize}
\item[(1)] $\mathcal{C}$ contains an expansion of odd type if and only if any minimal genus non-orientable spanning surface for $K$ is
boundary-incompressible.
\item[(2)] $\mathcal{C}$ contains no expansion of odd type if and only if any minimal genus non-orientable spanning surface for $K$ is
boundary-compressible.
\end{itemize}
\end{theorem}

The following is immediate from Theorem \ref{thm:characterization}.

\begin{corollary}\label{cor:compare}
A $2$-bridge knot never has two minimal genus non-orientable spanning surfaces such that
one is boundary-incompressible and the other is boundary-compressible.
\end{corollary}

In Section \ref{diagram}, we give an algorithm to visualize a minimal genus non-orientable spanning surface for $2$-bridge knots.

\begin{theorem}\label{thm:conway}
Any $2$-bridge knot $K$ has a Conway
diagram $D$ such that a minimal genus
non-orientable spanning surface for $K$ is obtained as a checker-board surface on $D$.
\end{theorem}

\begin{example}
We note that such diagrams are not unique and not shortest in general.
Let $K=S(15,4)$ as in Example \ref{ex:2g+1} with $\gamma(K)=3$.
Then the Conway diagrams $[2,1,5,-1,3]$
and $[4,1,1,1,4]$ representing $K$ respectively yield a desired surface as a checker-board surface.
It is interesting to confirm Theorem \ref{thm:bessho}(2) for these surfaces.
\end{example}

\section{Proof of Theorem \ref{thm:bessho}}\label{sec:bessho}

Let $K$ be a $2$-bridge knot with a minimal genus non-orientable spanning surface $F$.
Let $E(K)=S^3-\mathrm{Int}\,N(K)$ be its exterior.
Then $F\cap N(K)$ can be assumed to be a collar neighborhood of $\partial F$ in $F$, and hence
we will use the same notation $F$ for $F\cap E(K)$.

\begin{lemma}\label{incomp}
$F$ is incompressible in $E(K)$.
\end{lemma}

\begin{proof}
Assume not.  Let $D$ be a compressing disk for $F$.
Then $\partial D$ is an orientation-preserving loop on $F$.
Let $F'$ be the resulting surface from $F$ by compressing along $D$.
Then $\chi(F')=\chi(F)+2$.
If $F'$ is disconnected, then it consists of a closed orientable component $F_1$ and a non-orientable component $F_2$
with $\partial F_2\ne \emptyset$.
Since $\beta_1(F_1)+\beta_1(F_2)=\beta_1(F)$ and $\beta_1(F_1)>0$, we have
$\beta_1(F_2)<\beta_1(F)$. This contradicts the minimality of $\beta_1(F)$.
If $F'$ is connected and non-orientable, then $\beta_1(F')=\beta_1(F)-2$, a contradiction.
Hence $F'$ is connected and orientable.
This means that $F'$ is a Seifert surface for $K$.
Then adding a small half-twisted band to $F'$ gives a non-orientable spanning surface $R$ for $K$
with $\beta_1(R)=\beta_1(F')+1=\beta_1(F)-1$, a contradiction.
\end{proof}

\begin{proof}[Proof of Theorem \ref{thm:bessho}]
(1) Let $p/q=r+[b_1,b_2,\dots,b_n]$ be an expansion of odd type of some fraction $p/q$ for $K$.
We assume that the length $n$ is minimal among all expansions of odd type of all fractions for $K$.
The surface obtained by plumbing $n$ bands corresponding to this expansion in the usual way
gives a non-orientable spanning surface for $K$ with the first betti number $n$.
Thus $\gamma(K)\le n$.

The argument to show $n\le \gamma(K)$ is divided into two cases according to the boundary-incompressibility of
a minimal genus non-orientable spanning surface $F$.

First assume that $K$ has a minimal genus non-orientable spanning surface $F$ which is boundary-incompressible.
Then it is isotopic to one of the surfaces obtained by
plumbing $k$ bands corresponding to some expansion $s+[a_1,a_2,\dots,a_k]$ of some fraction for $K$ 
with $s\in \mathbb{Z}$ and $|a_i|\ge 2$ for each $i$ by \cite[Theorem 1(b)]{HT}.
Hence $\gamma(K)=k$.
Since $F$ is non-orientable, this expansion $s+[a_1,a_2,\dots,a_k]$ must be of odd type.
Thus $n\le k$ by the minimality of $n$, and hence we have $n\le \gamma(K)$.

Next, assume that any minimal genus non-orientable spanning surface $F$ for $K$ is boundary-compressible.
Let $D$ be a boundary-compressing disk such that $\partial D=\alpha\cup \beta$, where
$D\cap F=\alpha$ is a properly embedded essential arc in $F$ and $D\cap \partial E(K)=\beta$.
Then $\beta$ intersects $\partial F$ in two points.
If these two points have distinct signs (after orienting $\beta$ and $\partial F$ suitably),
then $\beta$ and a subarc of $\partial F$ bound a disk $\delta$ in $\partial E(K)$.
Thus $D\cup \delta$, pushed away from $\partial E(K)$ slightly,
gives a compressing disk for $F$, which contradicts Lemma \ref{incomp}.
Hence $\beta $ intersects $\partial F$ twice in the same direction.
Let $F_1$ be the surface obtained by boundary-compressing $F$ along $D$.
From the above observation, $F_1$ is a connected surface with connected boundary.
Also, we see $\beta_1(F_1)=\beta_1(F)-1$.

\begin{claim}
$F_1$ is incompressible in $E(K)$.
\end{claim}

\begin{proof}
Let $N(D)=D\times [-1,1]$ be a product neighborhood of $D$ such that
$N(D)\cap F=\partial N(D)\cap F=\alpha \times [-1,1]$.
Then $F_1=(F-N(D)\cap F)\cup (D\times \{-1,1\})$.
If $F_1$ is compressible, then it has a compressing disk $E$ disjoint from $N(D)$.
Since $F$ is incompressible, $\partial E$ bounds a disk $E'$ in $F$.
We can choose $E'$ disjoint from the disk $\alpha\times [-1,1]$.
Thus $\partial E$ bounds a disk in $F_1$, a contradiction. 
\end{proof}

If $F_1$ is orientable, then it is boundary-incompressible.
(Any orientable incompressible surface in $E(K)$ is boundary-incompressible if it has a connected boundary.)
If $F_1$ is non-orientable and boundary-compressible, we continue a boundary-compression.
Thus for some $\ell>0$ we have a sequence of incompressible surfaces
$F=F_0\to F_1 \to F_2 \to \dots \to F_\ell$
where $\beta_1(F_{i})=\beta_1(F_{i-1})-1$ for $i=1,2,\dots,\ell$ and $F_\ell$ is boundary-incompressible.
By \cite[Proposition 2]{HT}, $\partial F_\ell$ runs once longitudinally on $\partial E(K)$.
 
If $F_\ell$ is orientable, then it has minimal genus \cite[Corollary]{HT}.
Note that $\beta_1(F_\ell)=\beta_1(F)-\ell=\gamma(K)-\ell$ and that
$F_\ell$ corresponds to the unique expansion $r'+[b_1',b_2',\dots,b_m']$ with each $b_i'$ even.
Then $\beta_1(F_\ell)=m$.
Since $K$ admits an odd type expansion of length $m+1$,
$n\le m+1=\beta_1(F_\ell)+1=\gamma(K)-\ell+1\le \gamma(K)$.
Thus we have $n\le \gamma(K)$, and so $\gamma(K)=n$ and $\ell=1$.

If $F_\ell$ is non-orientable, then $F_\ell$ is isotopic to some surface
obtained by plumbing $k$ bands and $n\le k$ as before.
Since $\beta_1(F_\ell)=k$, $n\le k=\gamma(K)-\ell<\gamma(K)\le n$, a contradiction.
Thus such a case never occurs.

(2) If a minimal genus non-orientable spanning surface $F$ is boundary-compressible, then
the above argument shows that boundary-compressing $F$ gives a minimal genus Seifert surface $Q$ for $K$.
Then $F$ is obtained from $Q\subset E(K)$ by attaching a band $b=[0,1]\times [0,1]\subset \partial E(K)$ such that
$b\cap Q=[0,1]\times \{0,1\}$.
In fact, since $\partial F$ runs once longitudinally on $\partial E(K)$,
there are only two possibilities for $b$ as shown in Figure \ref{band}.
Thus $F$ is obtained from $Q$ by adding a M\"{o}bius band locally as desired.
\end{proof}

\section{Calculation by the modular diagram}\label{diagram}

\begin{figure}[ht]
\includegraphics*[scale=0.9]{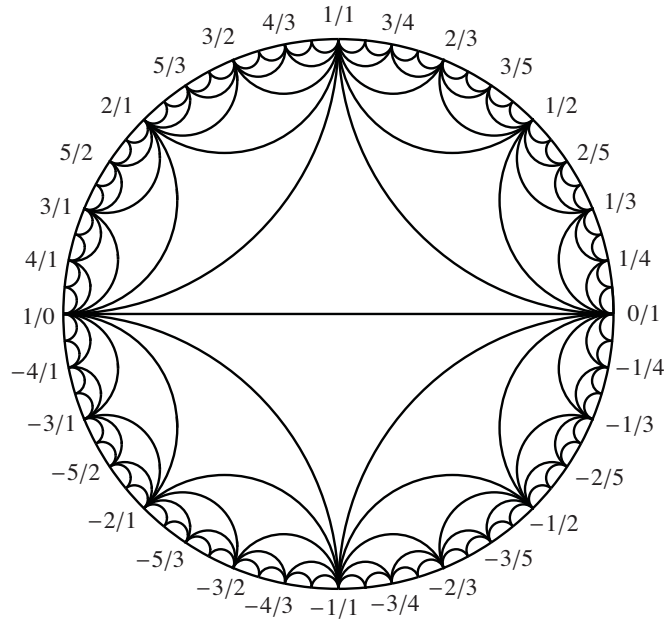}
\caption{The modular diagram $\mathcal{D}$}\label{ht}
\end{figure}

We use the modular diagram $\mathcal{D}$ as shown in Figure \ref{ht} to compute
the crosscap numbers of $2$-bridge knots.
This diagram comes from the action of $PSL(2,\mathbb{Z})$ on the hyperbolic plane.
(But $\mathcal{D}$ is distorted to space the vertices evenly along the circle.)

The vertices are labelled with $\mathbb{Q}\cup\{1/0\}$, inductively:
Start with $1/0$ and $0/1$ at the ends of the \lq horizontal\rq\! edge.
If two vertices of an triangle are already labelled with 
$a/b$ and $c/d$, then the remaining vertex of the triangle is 
labelled $(a+c)/(b+d)$. 
(This is the rule to label the vertices on the upper circle 
only, and for those on the lower circle, regard $1/0$ and 
$0/1$ as $-1/0$ and $-0/1$.)

We call the third vertex the \textit{child\/} of the first two 
vertices, which themselves are called the \textit{parents}, and call the edge connecting the parents the
\textit{longest side\/} of a triangle.
Note that two vertices $a/b$ and $c/d$ are connected by an edge if and only if $|ad-bc|=1$. 

We will identify a vertex with the corresponding label for convenience.
In fact, all rational numbers appear on the circle with the usual order.
That is, if the vertices $u$ and $v$ correspond to $a/b$ and $c/d$, respectively,
and if $a/b< c/d$ in $\mathbb{Q}$, then $u$ and $v$ lie on the circle with the counterclockwise orientation
in the order $u,v$.

\begin{figure}[ht]
\includegraphics*[scale=1]{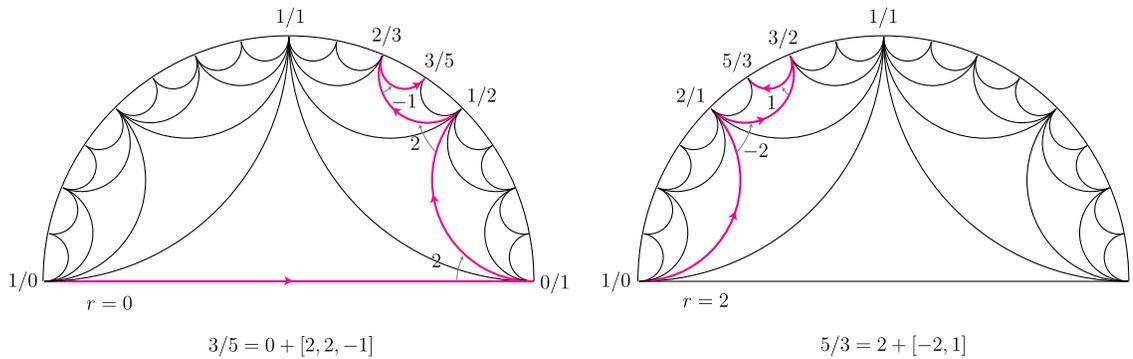}
\caption{Edge-paths from $1/0$ to $p/q$}\label{fig:path}
\end{figure}

An edge-path from $1/0$ to $p/q$ in $\mathcal{D}$ corresponds uniquely to an expansion
$p/q=r+[b_1,b_2,\dots,b_n]$, where
the partial sums $p_i/q_i=r+[b_1,b_2,\dots,b_i]$  ($p_0/q_0=r$)
are the successive vertices of the edge-path.
At the vertex $p_{i-1}/q_{i-1}$ the path turns left or right across $|b_i|$
triangles, left if $b_i>0$ and right if $b_i<0$.
See Figure \ref{fig:path}.
If an edge-path corresponds to an expansion of odd type, then the path is also said to be of odd type.

\begin{figure}[ht]
\includegraphics*[scale=0.35]{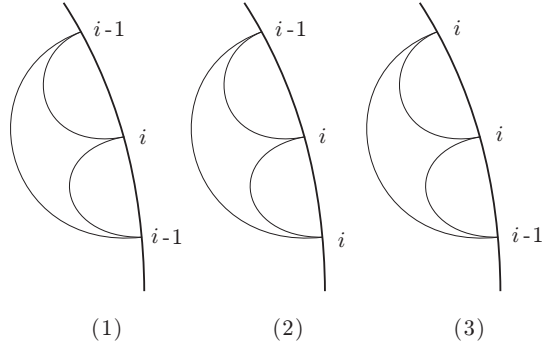}
\caption{Triangles}\label{triangle}
\end{figure}

We assign the \textit{depth\/} $d(v)$ to each vertex $v$ on $\mathcal{D}$.
First, set $d(v)=0$ for $v\in \mathbb{Z}\cup \{1/0\}$.
If a triangle in $\mathcal{D}$ has vertices $a/b$, $c/d$ and $(a+c)/(b+d)$,
then define the depth of the child from those of parents by setting
$$d((a+c)/(b+d))=\min\{d(a/b), d(c/d)\}+1.$$
Thus all vertices can be assigned the depths.
Notice that if two vertices $u$ and $v$ are connected by an edge in $\mathcal{D}$, then
$|d(u)-d(v)|\le 1$.
Also there are only three kinds of triangles in $\mathcal{D}$ as shown in Figure \ref{triangle},
where depths of vertices are indicated, except triangles with vertices $\{1/0,n,n+1\}$ of depth $0$ where $n\in \mathbb{Z}$.

\begin{lemma}\label{depth}
Let $\ell$ be the length of a shortest edge-path from $1/0$ to $p/q\ne 1/0$.
Then $d(p/q)=\ell-1$.
\end{lemma}

\begin{proof}
By the definition of depth, there is an edge-path from $1/0$ to $p/q$ of length $d(p/q)+1$.
Thus $\ell \le d(p/q)+1$.
Conversely, let $\xi$ be a shortest edge-path from $1/0$ to $p/q$.
Recall that  $|d(u)-d(v)|\le 1$ for any consecutive vertices $u, v$ on $\xi$.
In particular, $d(1/0)=d(v_0)=0$, where $v_0$ is the second vertex on $\xi$, and hence
$d(p/q)\le \ell-1$.  Thus we have $d(p/q)=\ell-1$.
\end{proof}



The next theorem gives a criterion for a shortest edge-path in terms of depths.

\begin{theorem}\label{criterion}
An edge-path $\xi: 1/0 \to v_0 \to v_1 \to v_2 \to \dots \to v_n=p/q\ne 1/0$
is shortest if and only if $d(v_i)=i$ for $i=0,1,2,\dots,n$.
\end{theorem}

\begin{proof}
Assume that $\xi$ is shortest.
Then $d(v_n)=n$ by Lemma \ref{depth}.
Since the depth can increase by at most one along $\xi$ and $d(v_0)=0$, $d(v_i)=i$ for each $i$. 

Conversely, since $d(v_n)=n$, any edge-path from $1/0$ to $p/q$ has length at least $n+1$ by Lemma \ref{depth}.
Thus we can conclude that $\xi$ is shortest.
\end{proof}


\begin{proof}[Proof of Theorem \ref{thm:algorithm}]
Let $p/q=r+[b_1,b_2,\dots,b_n]$ be an expansion fully reduced by the reductions in the statement of Theorem \ref{thm:algorithm}.
Let $\xi: 1/0 \to v_0 \to v_1 \to \dots \to v_n=p/q$
be the edge-path corresponding to the expansion.
Suppose that the expansion is not of minimal length, that is, $\xi$ is not shortest.
Then the sequence ${\mathcal S}$ of depths $d(v_0), d(v_1), d(v_2),\dots,d(v_n)$ is not strictly increasing
by Theorem \ref{criterion}.
Notice that $d(v_0)=0$.

First, suppose that ${\mathcal S}$ contains $i-1,i,i-1$ as a subsequence.
Then we can see that there is a triangle of Figure \ref{triangle}(1) such that $\xi$ runs along
the shorter two edges.
This means that some coefficient $b_j$ is $\pm 1$, a contradiction.

If ${\mathcal S}$ contains $0,0$, then $b_1=\pm 1$, a contradiction.
Thus ${\mathcal S}$ contains $i-1, i,i,\dots,i$ ($i$ is repeated $k\,(\ge 2)$ times) for some $i\ge 1$.
We choose $i$ minimal among such subsequences of ${\mathcal S}$.

Let $u_1$ and $u_2$ be the depth $i$ vertices on $\xi$, appearing in the order $u_2$, $u_1$.
We can suppose that the vertex before $u_2$ on $\xi$ has depth $i-1$.
There is the unique triangle $T_1$ which contains the edge between $u_1$ and $u_2$ as one of two shorter edges.
Without loss of generality, we can assume that $T_1$ has the form of Figure \ref{triangle}(3).
Let $w_1$ be the remaining vertex of $T_1$.
If $u_2$ is the child of $\{u_1,w_1\}$, then $\xi$ contains the edges $w_1\to u_2\to u_1$.
Then some $b_j=-1$, a contradiction.
Hence $u_1$ is the child of $\{u_2,w_1\}$.
Let $T_2$ be the (unique) triangle sharing the edge between $u_2$ and $w_1$ with $T_1$, and let $u_3$
be the remaining vertex of $T_2$.  See Figure \ref{232}(1).
(Since $d(w_1)=i-1$ and $d(u_2)=i$, $u_3$ is located in this position.)
Then $d(u_3)=i$ or $i-1$.
If $d(u_3)=i$, then $\xi$ contains $w_1\to u_2\to u_1$, so some coefficient is $-1$.
Thus $d(u_3)=i-1$.
If $\xi$ contains $w_1\to u_2\to u_1$, then some coefficient is $-1$ again.
Hence $\xi$ contains $u_3\to u_2\to u_1$.

\begin{figure}[ht]
\includegraphics*[scale=0.45]{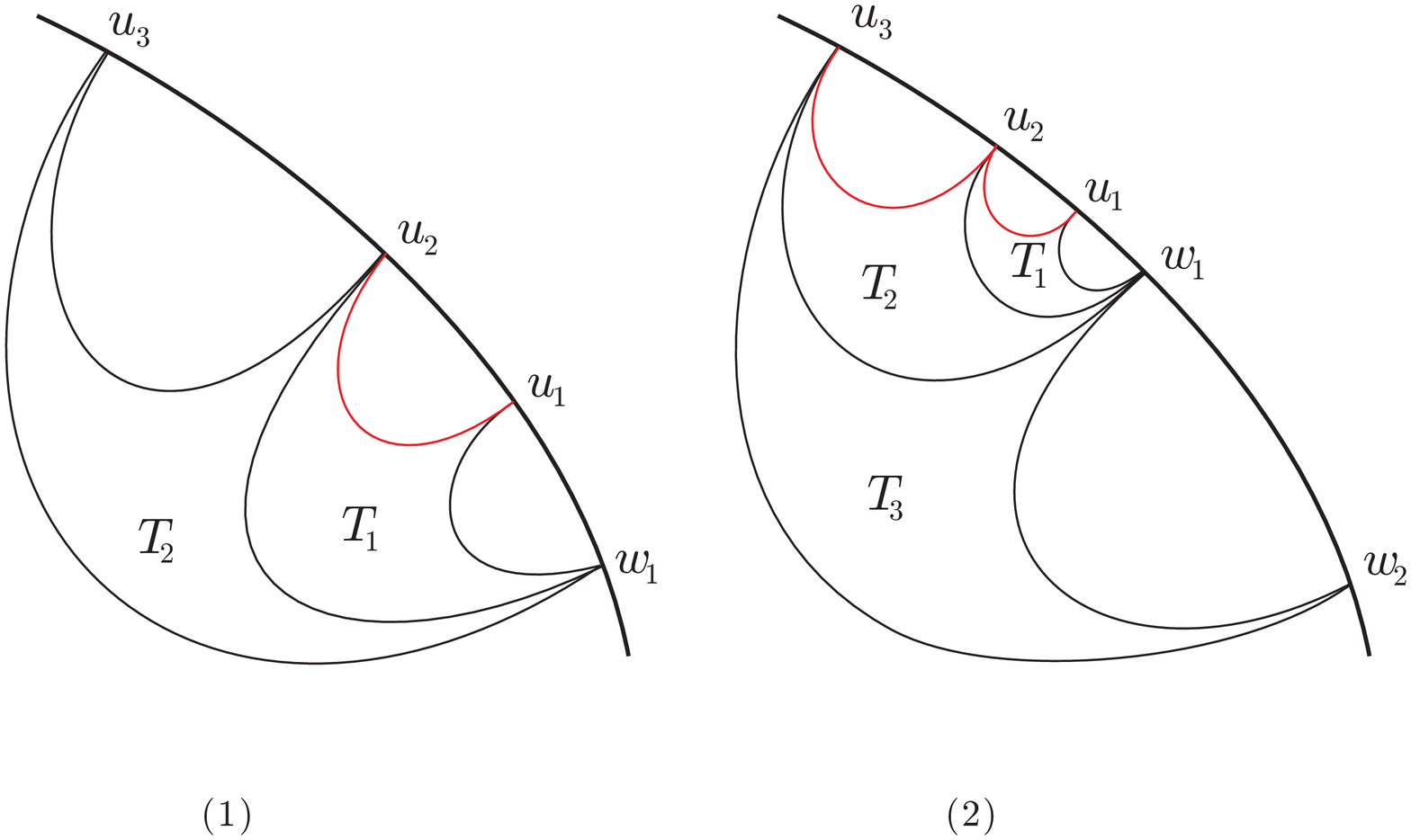}
\caption{Paths $u_2\to u_1$ and $u_3\to u_2\to u_1$}\label{232}
\end{figure}

If $i=1$, then $d(u_3)=d(w_1)=0$.
By the minimality of $i$,  $u_3=v_0$.
Thus $\xi$ contains $1/0\to u_3\to u_2\to u_1$, so $b_1=b_2=-2$, a contradiction.
(If $T_1$ has the form of Figure \ref{triangle}(2), then we encounter $1$ or $2,2$.)


Suppose $i\ge 2$.
Let $T_3$ be the triangle sharing the edge between $u_3$ and $w_1$ with $T_2$, and let $w_2$ be the remaining
vertex of $T_3$.
If $u_3$ is the child of $\{w_1,w_2\}$, then $d(w_2)=i-2$, and $\xi$ contains $w_2\to u_3\to u_2\to u_1$.
Then we have $-2,-2$ in the coefficients, a contradiction.
Hence $w_1$ is the child of $\{u_3,w_2\}$, and $d(w_2)=i-2$.  See Figure \ref{232}(2).
Let $T_4$ be the triangle sharing the edge between $u_3$ and $w_2$ with $T_3$, and let $u_4$ be the remaining vertex of $T_4$.
Then $u_{3}$ is the child of $\{u_{4},w_2\}$, and so $d(u_{4})=i-1$ or $i-2$.
If $d(u_{4})=i-1$, then $\xi$ contains $w_2\to u_3 \to u_2 \to u_1$.
Then we have $-2,-2$ in the coefficients.
Hence $d(u_{4})=i-2$.
If $\xi$ contains $w_2\to u_3\to u_2\to u_1$, then
we have $-2,-2$ again.
Thus $\xi$ contains $u_{4}\to u_3\to u_2\to u_1$.
See Figure \ref{2332}(1).
If $i=2$, then $d(u_4)=d(w_2)=0$ and hence $u_4=v_0\,(\ne 1/0)$.
Thus $\xi$ contains $1/0\to u_4\to u_3\to u_2\to u_1$, so
$b_1=-2$, $b_2=-3$ and $b_3=-2$, a contradiction.
Thus $i\ge 3$.

\begin{figure}[ht]
\includegraphics*[scale=0.5]{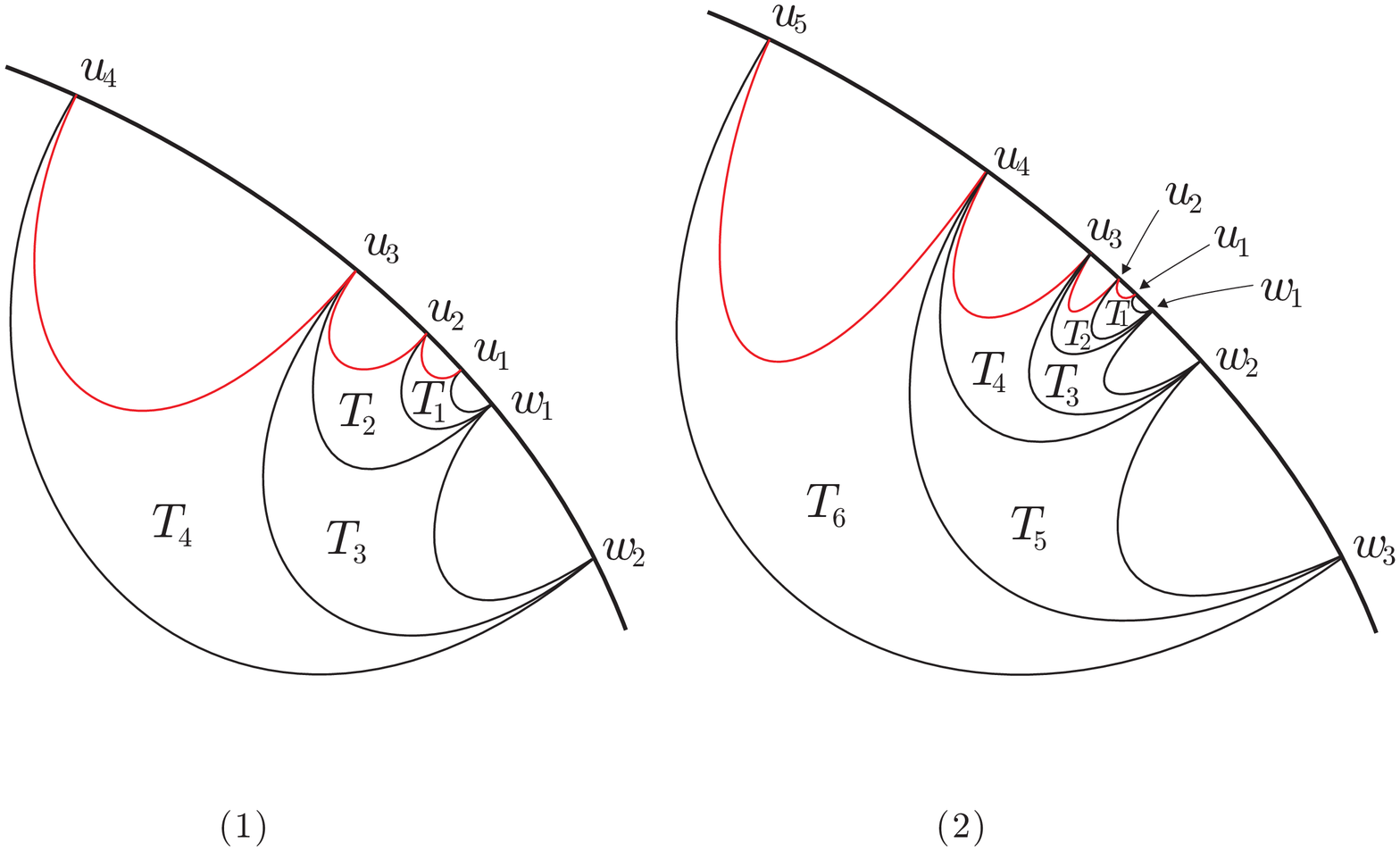}
\caption{Paths $u_{4}\to u_3\to u_2\to u_1$ and $u_5\to u_{4}\to u_3\to u_2\to u_1$}\label{2332}
\end{figure}

Continuing this process, we obtain the triangles $T_1,T_2,\dots,T_{2i}$, where
$T_{2m-1}$ is the form of Figure \ref{triangle}(3), and $T_{2m}$ is of Figure \ref{triangle}(1). 
Also, $T_1$ has vertices $\{u_1,u_2,w_1\}$, $T_{2m-1}$ ($m\ge 2$) has vertices $\{w_{m-1},w_{m},u_{m+1}\}$, and
$T_{2m}$ $(m\ge 1)$ has vertices $\{u_{m+1},u_{m+2},w_{m}\}$, and $d(u_1)=i$, $d(u_j)=i-j+2$ for $j\ge 2$, $d(w_j)=i-j$ for $j\ge 1$.
The path $\xi$ contains the edges $u_{i+2}\to u_{i+1}\to \dots \to u_2\to u_1$.
Since $d(u_{i+2})=0$, $u_{i+2}=v_0$.
Then $-2,-3,\dots,-3,-2$ appears in the coefficients, a contradiction.
(If $T_1$ has the form of Figure \ref{triangle}(2), then we encounter $1$ or $2,3,3,\dots,3,2$ in the coefficients.)
\end{proof}

\begin{lemma}\label{lem:key}
Let $p/q$ and $p'/q'$ be two fractions for a $2$-bridge knot $K$.
Then there exists a one-to-one correspondence between the set of all expansions
for $p/q$ and that of $p'/q'$, such that the correspondence preserves both length and type of expansions.
\end{lemma}


\begin{proof}
Since both $p/q$ and $p'/q'$ represent the same knot, $q=q'$ and (i) $p\equiv p' \pmod{q}$ or (ii) $p p'\equiv 1\pmod{q}$.
Suppose $p/q=p'/q +s$, where $s\in {\mathbb Z}$.
Then for any expansion $r+[a_1,a_2, \ldots, a_n]$ of $p/q$, we can associate $r-s+[a_1,a_2, \ldots, a_n]$ of $p'/q$.
Therefore, it suffices to establish a one-to-one correspondence only for $p/q$ and $p'/q$ lying between $0$ and $1$.
Under such a restriction, suppose that $p p'\equiv 1\pmod{q}$.
Then it is well known that
if $[a_1,a_2,\ldots, a_n]=p/q$ then $[a_n, \ldots, a_2,a_1]= p'/q$.
Thus we can establish a required one-to-one correspondence.
\end{proof}

\begin{proof}[Proof of Theorem \ref{thm:main}]
Let $p/q$ be a fraction corresponding to $K$.
Let $n$ be the minimal length of an expansion of $p/q$ obtained by the reduction of Theorem \ref{thm:algorithm}.
If another fraction for $K$ admits an expansion of length shorter than $n$, then
$p/q$ has an expansion of the same length by Lemma \ref{lem:key}.
This contradicts the minimality of $n$.
\end{proof}

Next, we consider when there exists a shortest edge-path from $1/0$ to $p/q$, which is of odd type.
A \textit{rectangle\/} in $\mathcal{D}$ is the union of two triangles sharing one edge.

\begin{lemma}\label{rectangle}
Assume that a shortest edge-path
$\xi$ from $1/0$ to $p/q$ is not of odd type.
If there is a rectangle in $\mathcal{D}$ containing
two successive edges of $\xi$, then there is a shortest edge-path of odd type from $1/0$ to $p/q$.
\end{lemma}

\begin{proof}
Without loss of generality, we can assume that
the successive two edges on $\xi$, $v_{i-2}\to v_{i-1}$ and $v_{i-1}\to v_i$
lie on a rectangle which is the union of two triangles whose vertices $\{v_{i-2},v_{i-1},v_{i-1}'\}$ and
$\{v_{i-1},v_{i-1}',v_i\}$.
Then we replace these two edges by $v_{i-2}\to v_{i-1}'$ and $v_{i-1}'\to v_i$.
This will change $\xi$ to a new shortest edge-path of odd type.
For, if $v_{i-2}\ne 1/0$, then 
the expansion corresponding to $\xi$ changes
from $r+[\dots,a,\pm 2,b,\dots]$ to $r+[\dots,a\mp 1,\mp 2,b\mp 1,\dots]$, or $r+[\dots,a,\pm 2]$ to
$r+[\dots,a\mp 1,\mp 2]$.
If $v_{i-2}=1/0$, then $v_{i-1}\in \mathbb{Z}$.
Thus $\xi$ corresponds to an expansion
$p/q=v_{i-1}+[\pm 2,a,\dots]$.
Then the new edge-path corresponds to $v_{i-1}\pm 1+[\mp 2,a\mp 1,\dots]$.
\end{proof}

The deformation used in the proof of Lemma \ref{rectangle} is referred to as
a \textit{rectangle move}.

\begin{proof}[Proof of Theorem \ref{thm:shortest}]
Let 
$$\xi: 1/0\to v_0 \to v_1 \to v_2 \to \dots \to v_n=p/q,$$
$$\xi':  1/0\to v_0' \to v_1' \to v_2' \to \dots \to v_n'=p/q$$
be two shortest edge-paths from $1/0$ to $p/q$.
Recall that each of $v_i$ and $v_i'$ has depth $i$ for any $i$ by Theorem \ref{criterion}.
In particular, $v_0,v_0'\in \mathbb{Z}$.

Suppose $v_0\ne v_0'$.
We may assume that $v_0<v_0'$.
Since $p/q\in (v_0-1,v_0+1)$ and $p/q\in (v_0'-1,v_0'+1)$,
$v_0'-v_0=1$.
In fact, $v_1$ and $v_1'$ lie in the interval $(v_0, v_0')$.
Let $u$ be the child of $\{v_0,v_0'\}$.
Then $d(u)=1$.
If neither of $v_1$ nor $v_1'$ is $u$, then we would have $p/q<u<p/q$, a contradiction.
Hence we may assume $v_1'=u$.
Then a single application of rectangle move on $\xi'$ changes the edges
$1/0\to v_0'\to v_1'=u$ to $1/0\to v_0\to u$.

Suppose that $v_i=v_i'$ for $0\le i\le k$ and $v_{k+1}\ne v_{k+1}'$ for some $k\ge 0$.
If $v_{k+1}<v_k=v_k'<v_{k+1}'$ or $v_{k+1}> v_k=v_k'>v_{k+1}'$, then
we would have $p/q<v_k=v_k'<p/q$, a contradiction. 
Hence we consider only the case where $v_{k+1},v_{k+1}'>v_k=v_k'$.
The case where $v_{k+1},v_{k+1}'<v_k=v_k'$ is similar.
Without loss of generality, assume $v_{k+1}'<v_{k+1}$.
Then $p/q$ lies in the interval $(v_{k+1}',v_{k+1})$.
If there is a vertex $u$ with depth $k+1$ inside the interval, then
we have $p/q<u<p/q$, a contradiction.
Hence there are a triangle $\Delta_1$ whose vertices
are $v_{k},v_{k+1}$ and $v_{k+1}'$, and a triangle $\Delta_2$
whose vertices are $v_{k+1},v_{k+1}'$ and $w$.
Here, $v_{k+1},v_{k+1}'$ are the parents of $w$.
If $w=p/q$, then $v_{k+2}=v_{k+2}'$.
Then $\xi$ can be changed to $\xi'$ by the rectangle move on $\Delta_1\cup \Delta_2$.
Otherwise $p/q$ lies in $(v_{k+1}',w)$ or $(w,v_{k+1})$.
Then $w=v_{k+2}$ in the former case, and $w=v_{k+2}'$ in the latter.
After the rectangle move on $\Delta_1\cup \Delta_2$ to $\xi$ or $\xi'$,
we have $v_{k+1}=v_{k+1}'$.
Thus, $\xi$ can be changed to $\xi'$ gradually.
\end{proof}

\begin{proof}[Proof of Theorem \ref{thm:final}]
If a shortest expansion of $p/q$ obtained by Theorem \ref{thm:algorithm} contains an odd
coefficient, then $\gamma(K)=n$ by Theorems \ref{thm:bessho} and \ref{thm:main}.
Let $\xi$ be the corresponding edge-path and assume that $\xi$ is of even type.
If the expansion contains a coefficient $\pm 2$, then there is a rectangle in $\mathcal{D}$ containing two successive edges of $\xi$.
Hence a rectangle move creates another shortest edge-path of odd type by Lemma \ref{rectangle}.
Then $\gamma(K)=n$ as above.
Otherwise $\xi$ is the unique shortest edge-path from $1/0$ to $p/q$ by Theorem \ref{thm:shortest}.

If another fraction for $K$ admits an expansion of odd type of length $n$, then
$p/q$ also admits an expansion of odd type of length $n$ by Lemma \ref{lem:key}.
Thus there is no expansion of odd type of length $n$, and so $\gamma(K)=n+1$.
\end{proof}

\begin{proof}[Proof of Corollary \ref{thm:2g+1}]
It is well known that even type expansions for a $2$-bridge knot
are unique modulo integer parts, and the length equals $2g(K)$.
If the expansion does not contain $\pm 2$, then it is shortest by Theorem \ref{thm:algorithm}.
Hence $\gamma(K)=2g(K)+1$ by Theorem \ref{thm:final}.
If the expansion contains $\pm 2$, then there is another expansion of odd type
with the same length by Lemma \ref{rectangle}.
This means that $\gamma(K)\le 2g(K)$.
\end{proof}

\begin{proof}[Proof of Theorem \ref{thm:characterization}]
Let $F$ be a minimal genus non-orientable spanning surface for $K$.
By Lemma \ref{incomp}, $F$ is incompressible.

First, assume that the minimal length $n$ of expansions of all fractions for $K$ is
realized by an expansion of odd type.

If $F$ is boundary-compressible, then
boundary-compression yields a minimal genus Seifert surface $S$ for $K$ with $\beta_1(S)=\beta_1(F)-1$
as in the proof of Theorem \ref{thm:bessho}.
Note that $S$ is isotopic to a plumbing surface corresponding to the unique expansion with only even coefficients.
In particular, such expansion has length $n-1$.
This contradicts the minimality of $n$.
Hence $F$ is boundary-incompressible.
(In this case $F$ is isotopic to a plumbing surface by \cite{HT}.)

Next, assume that only expansions of even type realize the minimal length $n$.
If $F$ is boundary-incompressible, then $F$ is isotopic to a plumbing surface which corresponds to some
expansion of odd type, which has length $k\ge n+1$.
Indeed, $k=n+1$ by the minimality of $\beta_1(F)$.
If $F$ corresponds to an expansion $r+[b_1,b_2,\dots,b_{n+1}]$, then $|b_i|\ge 2$ for each $i$ by \cite{HT}.
This expansion is not shortest, and hence the sequence $b_1,b_2,\dots,b_{n+1}$ contains
a subsequence $2\varepsilon,3\varepsilon,\dots,3\varepsilon,2\varepsilon$, $\varepsilon=\pm 1$, where
the number of $3\varepsilon$ may be zero, by Theorem \ref{thm:algorithm}.
But such expansion can be reduced as shown in Theorem \ref{thm:algorithm}.
In particular, we have a shorter expansion of odd type, a contradiction.
Thus $F$ must be boundary-compressible.
\end{proof}

\begin{proof}[Proof of Theorem \ref{thm:conway}]
Let $K=S(q,p)$ with $\gamma(K)=\gamma$.
We omit the trivial case $\gamma=1$.
Let $C=[a_1, b_1, a_2, b_2,\cdots]$ be a shortest expansion for $p/q$ among those of odd type.
To be precise take $[a_1, b_1, \cdots, a_{(\gamma+1)/2}]$ if $\gamma$ is odd and if otherwise, take $[a_1, b_1, \cdots, b_{\gamma/2}]$.
Figure \ref{fig:check}(1) represents a Conway diagram for $K=S(q,p)$ of length $\gamma$.
Deform it to Figure \ref{fig:check}(3) through (2)
corresponding to the expansion $C'=[a_1-1, -1, b_1, 1, a_2, -1, b_2, 1, \cdots]$,
which ends with $a_{(\gamma+1)/2}+1$ (resp.\ $1$) if $\gamma$ is odd (resp.\ even). 
\begin{figure}[ht]
\includegraphics*[scale=0.8]{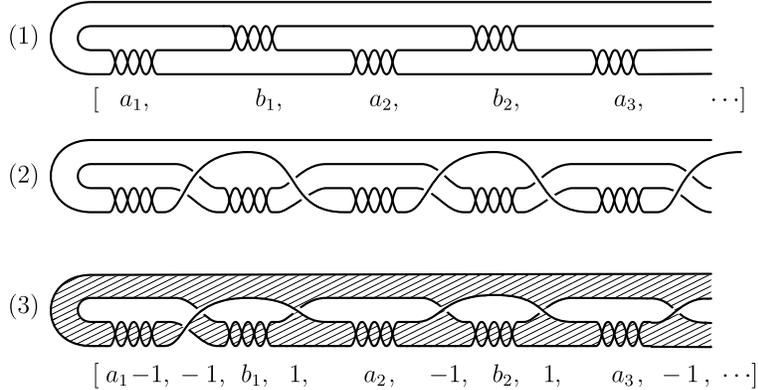}
\caption{Deformation of Conway diagram}\label{fig:check}
\end{figure}
Then we have a desired Conway diagram with a checkerboard  surface $F$ with $\beta(F)=\gamma$.
We  remark it is not preferable if $a_1-1=0$ or $a_{(\gamma+1)/2}+1=0$.
 In that case, apply the following:
Take the mirror image of $K$ and thus change all the signs of $C$, apply the deformation in Figure \ref{fig:check},
and the take the mirror image again.
Then we obtain a desired Conway diagram for $K$.
(This modification works since $C$ has at most one coefficient $\pm 1$ because of the minimality of $\gamma$.)
\end{proof}

\section{Examples}\label{sec:example}

In \cite{U}, it is proved that any positive integer can appear as the crosscap number of
some pretzel knot. 
We can show that such examples can be found among $2$-bridge knots.

As seen in Section \ref{sec:table}, among $2$-bridge knots up to $12$ crossings, exactly $7_4$, $8_3$,
$9_5$, $10_3$, $11a_{343}$, $11a_{363}$, $12a_{1166}$ and $12a_{1287}$ do not have a shortest expansion of odd type and hence satisfy the
equality $\gamma(K)=2g(K)+1$.
The next example also gives an infinite series of such $2$-bridge knots, as a generalization of \cite{MY}.

\begin{example}
Let $K_{m,n}$ be the $2$-bridge knot corresponding to $[m,4,4,\dots,4]$ of length $n$
for any $n\ge 1$. When $n=1$, $K_{m,1}=S(m,1)$.

If $m$ is odd and $m\ge 3$, then this expansion is shortest by Theorem \ref{thm:algorithm}.
Thus $\gamma(K_{m,n})=n$ by Theorem \ref{thm:final}.
Also, if $m\ne m'$, then $K_{m,n}$ and $K_{m',n}$ have distinct denominators, and hence
they are not equivalent.
Thus we have infinitely many $2$-bridge knots $K_{m,n}$ with $\gamma(K_{m,n})=n$ for any $n\ge 1$.

If $m$ and $n$ are even and $m\ge 4$, $n\ge 2$, then
$g(K_{m,n})=n/2$.
By Corollary \ref{thm:2g+1}, $\gamma(K_{m,n})=n+1$.
Also, distinct $m$'s give distinct knots.
Thus we have infinitely many $2$-bridge knots $K_{m,n}$ satisfying the equality $\gamma(K_{m,n})=2g(K_{m,n})+1$.
\end{example}

The Murasugi sum of two minimal genus Seifert surfaces gives 
a minimal genus Seifert surface \cite{G}.
Finally, we give the examples of $2$-bridge knots, showing that
an analogous statement does not hold in non-orientable case.
This is a generalization of Bessho's example \cite{B}.

\begin{example}
For any odd integer $m\ge 3$, let $K_m$ be the $2$-bridge knot corresponding to $[m,2]$.
Then $\gamma(K)=2$ and its minimal genus non-orientable spanning surface $F$ is
obtained by plumbing two bands with $m$ and $2$ half-twists respectively.
Let $R$ be the Murasugi sum of two copies $F_1$ and $F_2$ of $F$.
Here, the plumbing disks are chosen to lie in the band with $2$ half-twists of $F_1$ and 
in the band with $2$ half-twists of $F_2$.
Then $\beta_1(R)=4$.
But $\partial R$ is the $2$-bridge knot corresponding to $[m,2,2,m]$.
Since $[m,2,2,m]=[m-1,-3,m-1]$, it has crosscap number $3$.
Also, distinct $m$'s give inequivalent knots.
Thus the Murasugi sum of two minimal genus non-orientable spanning surfaces is not necessarily minimal genus.
\end{example}

\section{Table}\label{sec:table}

Here is the table of crosscap numbers of $362$ $2$-bridge knots up to
$12$ crossings.
The numbering of knots with $10$ or less crossings follows that of \cite{Ro}.
For $11$, $12$ crossings knots, we have used Dowker-Thistlethwaite notation.
The last column gives a minimal length subtractive continued fraction expansion of $p/q$.
We chose them to be of odd type except for the ones (indicated by $\ast$) where 
the shortest expansion is unique and of even type.
We referred to \cite{BZ} for $2$-bridge knots up to $10$ crossings, to 
\cite{CL} for those of $11$ and $12$ crossings, for which we 
also used a table compiled by David De Wit \cite{D}.

\bigskip
\begin{center}
\begin{tabular}{llllllll}
knot & $p/q$ & $\gamma$ & continued fraction & knot & $p/q$ & $\gamma$ & continued fraction \\
\hline
$3_1$ & $1/3$ & $1$ & $[3]$           & $4_1$ & $2/5$ & $2$ & $[3,2]$ \\
$5_1$ & $1/5$ & $1$ & $[5]$           & $5_2$ & $3/7$ & $2$ & $[2,-3]$ \\
$6_1$ & $2/9$ & $2$ & $[5,2]$     & $6_2$ & $4/11$ & $2$ & $[3,4]$ \\
$6_3$ & $5/13$ & $3$ & $[3,2,-2]$     & $7_1$ & $1/7$ & $1$ & $[7]$ \\
$7_2$ & $5/11$ & $2$ & $[2,-5]$       & $7_3$ & $4/13$ & $2$ & $[3,-4]$ \\
$7_4$ & $4/15$ & $3$ & $[4,4]^*$     & $7_5$ & $7/17$ & $3$ & $[2,-2,3]$ \\
$7_6$ & $7/19$ & $3$ & $[3,3,-2]$     & $7_7$ & $8/21$ & $3$ & $[3,3,3]$ \\
$8_1$ & $2/13$ & $2$ & $[7,2]$    & $8_2$ & $6/17$ & $2$ & $[3,6]$ \\
$8_3$ & $4/17$ & $3$ & $[4,-4]^*$     & $8_4$ & $5/19$ & $2$ & $[4,5]$ \\
$8_6$ & $10/23$ & $3$ & $[2,-3,3]$    & $8_7$ & $9/23$ & $3$ &  $[3,2,-4]$ \\
$8_8$ & $9/25$ & $3$ &  $[3,4,-2]$    & $8_9$ & $7/25$ & $3$ &  $[4,2,-3]$ \\
$8_{11}$ & $10/27$ & $3$ & $[3,3,-3]$ & $8_{12}$ & $12/29$ & $4$ & $[3,2,4,2]$ \\
$8_{13}$ & $11/29$ & $3$ & $[3,3,4]$  & $8_{14}$ & $12/31$ & $4$ & $[3,2,-2,2]$ \\
$9_1$ & $1/9$ & $1$ & $[9]$           & $9_2$ & $7/15$ & $2$ & $[2,-7]$ \\
$9_3$ & $6/19$ & $2$ & $[3,-6]$       & $9_4$ & $5/21$ & $2$ & $[4,-5]$ \\
$9_5$ & $6/23$ & $3$ & $[4,6]^*$     & $9_6$ & $5/27$ & $3$ & $[5,-2,2]$ \\
$9_7$ & $13/29$ & $3$ & $[2,-4,3]$    & $9_8$ & $11/31$ & $3$ & $[3,5,-2]$ \\
$9_9$ & $9/31$ & $3$ & $[3,-2,4]$     & $9_{10}$ & $10/33$ & $3$ & $[3,-3,3]$ \\
$9_{11}$ & $14/33$ & $3$ & $[2,-3,-5]$ & $9_{12}$ & $13/35$ & $3$ & $[3,3,-4]$ \\
$9_{13}$ & $10/37$ & $3$ & $[4,3,-3]$ & $9_{14}$ & $14/37$ & $3$ & $[3,3,5]$ \\
$9_{15}$ & $16/39$ & $4$ & $[2,-2,3,-2]$ & $9_{17}$ & $14/39$ & $3$ & $[3,5,3]$ \\
$9_{18}$ & $17/41$ & $4$ & $[2,-2,2,-3]$ & $9_{19}$ & $16/41$ & $4$ & $[3,2,-3,2]$ \\
$9_{20}$ & $15/41$ & $3$ & $[3,4,4]$ & $9_{21}$ & $18/43$ & $4$ & $[2,-3,-2,3]$ \\
$9_{23}$ & $19/45$ & $4$ & $[2,-3,-3,2]$ & $9_{26}$ & $18/47$ & $4$ & $[3,3,2,-3]$ \\
$9_{27}$ & $19/49$ & $4$ & $[3,2,-3,-3]$ & $9_{31}$ & $21/55$ & $4$ & $[3,3,3,3]$ \\
$10_1$ & $2/17$ & $2$ & $[9,2]$  & $10_2$ & $8/23$ & $2$ & $[3,8]$ \\
$10_3$ & $6/25$ & $3$ & $[4,-6]^*$   & $10_4$ & $7/27$ & $2$ & $[4,7]$ \\
$10_5$ & $13/33$ & $3$ & $[3,2,-6]$  & $10_6$ & $16/37$ & $3$ & $[2,-3,5]$ \\
$10_7$ & $16/43$ & $3$ & $[3,3,-5]$  & $10_8$ & $6/29$ & $2$ & $[5,6]$ \\
$10_9$ & $11/39$ & $3$ & $[4,2,-5]$  & $10_{10}$ & $17/45$ & $3$ & $[3,3,6]$ \\
$10_{11}$ & $13/43$ & $3$ & $[3,-3,4]$   & $10_{12}$ & $17/47$ & $3$ & $[3,4,-4]$ \\
$10_{13}$ & $22/53$ & $4$ & $[3,2,3,-4]$ & $10_{14}$ & $22/57$ & $4$ & $[3,2,-2,4]$ \\
\end{tabular}
\end{center}
\begin{center}
\begin{tabular}{llllllll}
knot & $p/q$ & $\gamma$ & continued fraction & knot & $p/q$ & $\gamma$ & continued fraction \\
\hline
$10_{15}$ & $19/43$ & $3$ & $[2,-4,-5]$  & $10_{16}$ & $14/47$ & $3$ & $[3,-3,-5]$ \\
$10_{17}$ & $9/41$ & $3$ & $[5,2,-4]$    & $10_{18}$ & $23/55$ & $4$ & $[2,-3,-2,4]$ \\
$10_{19}$ & $14/51$ & $3$ & $[4,3,5]$    & $10_{20}$ & $16/35$ & $3$ & $[2,-5,3]$ \\
$10_{21}$ & $16/45$ & $3$ & $[3,5,-3]$   & $10_{22}$ & $13/49$ & $3$ & $[4,4,-3]$ \\
$10_{23}$ & $23/59$ & $4$ & $[3,2,-3,3]$ & $10_{24}$ & $24/55$ & $4$ & $[2,-3,2,-3]$ \\
$10_{25}$ & $24/65$ & $4$ & $[3,3,-2,3]$ & $10_{26}$ & $17/61$ & $4$ & $[4,2,-2,3]$ \\
$10_{27}$ & $27/71$ & $4$ & $[3,3,3,-3]$  & $10_{28}$ & $19/53$ & $3$ & $[3,5,4]$ \\
$10_{29}$ & $26/63$ & $4$ & $[2,-2,3,4]$  & $10_{30}$ & $26/67$ & $4$ & $[3,2,-3,-4]$ \\
$10_{31}$ & $25/57$ & $4$ & $[2,-4,-2,3]$ & $10_{32}$ & $29/69$ & $4$ & $[2,-3,-3,-4]$ \\
$10_{33}$ & $18/65$ & $4$ & $[4,3,2,-3]$  & $10_{34}$ & $13/37$ & $3$ & $[3,6,-2]$ \\
$10_{35}$ & $20/49$ & $4$ & $[3,2,6,2]$   & $10_{36}$ & $20/51$ & $4$ & $[3,2,-4,2]$ \\
$10_{37}$ & $23/53$ & $4$ & $[2,-3,3,-2]$ & $10_{38}$ & $25/59$ & $4$ & $[2,-3,-4,2]$ \\
$10_{39}$ & $22/61$ & $4$ & $[3,4,-2,2]$  & $10_{40}$ & $29/75$ & $5$ & $[3,2,-2,2,-2]$ \\
$10_{41}$ & $26/71$ & $4$ & $[3,4,3,-2]$  & $10_{42}$ & $31/81$ & $5$ & $[3,3,2,-2,2]$ \\
$10_{43}$ & $27/73$ & $4$ & $[3,3,-3,-3]$ & $10_{44}$ & $30/79$ & $4$ & $[3,3,4,3]$ \\
$10_{45}$ & $34/89$ & $5$ & $[3,3,3,2,-2]$ & $11a_{13}$ & $28/61$ & $4$ & $[2,-6,-3,-2]$ \\
$11a_{59}$ & $20/43$ & $3$ & $[2,-7,-3]$ & $11a_{65}$ & $27/59$ & $4$ & $[2,-5,2,-2]$ \\
$11a_{75}$ & $36/83$ & $4$ & $[2,-3,4,3]$ & $11a_{77}$ & $55/131$ & $5$ & $[2,-3,-3,-3,-3]$ \\
$11a_{84}$ & $44/101$ & $5$ & $[2,-3,3,2,-2]$ & $11a_{85}$ & $47/107$ & $5$ & $[2,-4,-2,2,3]$ \\
$11a_{89}$ & $44/119$ & $5$ & $[3,3,-3,-2,2]$ & $11a_{90}$ & $23/87$ & $4$ & $[4,4,-2,-3]$ \\
$11a_{91}$ & $50/129$ & $5$ & $[3,2,-3,-3,-3]$ & $11a_{93}$ & $41/93$ & $4$ & $[2,-4,-5,-3]$ \\
$11a_{95}$ & $33/73$ & $4$ & $[2,-5,-3,2]$ & $11a_{96}$ & $50/121$ & $5$ & $[2,-2,3,3,3]$ \\
$11a_{98}$ & $18/77$ & $4$ & $[4,-4,-3,-2]$ & $11a_{110}$ & $35/97$ & $4$ & $[3,4,-3,-3]$ \\
$11a_{111}$ & $37/103$ & $4$ & $[3,5,3,3]$ & $11a_{117}$ & $49/117$ & $5$ & $[2,-3,-2,3,3]$ \\
$11a_{119}$ & $34/77$ & $4$ & $[2,-4,-5,-2]$ & $11a_{120}$ & $45/109$ & $5$ & $[2,-2,3,3,-2]$ \\
$11a_{121}$ & $50/119$ & $5$ & $[2,-3,-3,-3,2]$ & $11a_{140}$ & $17/65$ & $3$ & $[4,6,3]$ \\
$11a_{144}$ & $17/73$ & $4$ & $[4,-3,2,-2]$ & $11a_{145}$ & $22/83$ & $4$ & $[4,4,-3,-2]$ \\
$11a_{154}$ & $30/67$ & $4$ & $[2,-4,3,-2]$ & $11a_{159}$ & $46/111$ & $5$ & $[2,-2,2,-3,-3]$ \\
$11a_{166}$ & $14/59$ & $3$ & $[4,-5,-3]$ & $11a_{174}$ & $28/79$ & $4$ & $[3,5,-2,-3]$ \\
$11a_{175}$ & $41/105$ & $5$ & $[3,2,-4,-2,2]$ & $11a_{176}$ & $31/111$ & $5$ & $[4,2,-3,-2,2]$ \\
$11a_{177}$ & $21/97$ & $4$ & $[5,3,3,3]$ & $11a_{178}$ & $34/123$ & $5$ & $[4,2,-2,-3,-3]$ \\
$11a_{179}$ & $20/57$ & $3$ & $[3,7,3]$ & $11a_{180}$ & $25/89$ & $4$ & $[4,2,-4,-3]$ \\
$11a_{182}$ & $13/73$ & $4$ & $[6,2,-2,-3]$ & $11a_{183}$ & $34/115$ & $5$ & $[3,-3,-2,2,3]$ \\
$11a_{184}$ & $19/87$ & $4$ & $[5,2,-3,-3]$ & $11a_{185}$ & $30/109$ & $4$ & $[4,3,4,3]$ \\
$11a_{186}$ & $39/95$ & $5$ & $[2,-2,3,-2,2]$ & $11a_{188}$ & $14/67$ & $3$ & $[5,5,3]$ \\
$11a_{190}$ & $18/85$ & $4$ & $[5,3,-2,-3]$ & $11a_{191}$ & $19/83$ & $4$ & $[4,-3,-3,2]$ \\
$11a_{192}$ & $26/97$ & $4$ & $[4,4,3,-2]$ & $11a_{193}$ & $29/95$ & $4$ & $[3,-4,-3,-3]$ \\
$11a_{195}$ & $8/53$ & $3$ & $[7,3,3]$ & $11a_{203}$ & $11/63$ & $3$ & $[6,4,3]$ \\
$11a_{204}$ & $30/101$ & $4$ & $[3,-3,-4,-3]$ & $11a_{205}$ & $25/91$ & $4$ & $[4,3,-4,-2]$ \\
$11a_{206}$ & $7/47$ & $3$ & $[7,3,-2]$ & $11a_{207}$ & $26/85$ & $4$ & $[3,-4,-3,2]$ \\
$11a_{208}$ & $31/105$ & $5$ & $[3,-3,-2,2,-2]$ & $11a_{210}$ & $16/73$ & $4$ & $[5,2,-3,2]$ \\
$11a_{211}$ & $12/67$ & $4$ & $[6,2,-3,-2]$ & $11a_{220}$ & $23/85$ & $4$ & $[4,3,-3,2]$ \\
$11a_{224}$ & $27/89$ & $4$ & $[3,-3,3,3]$ & $11a_{225}$ & $11/53$ & $3$ & $[5,5,-2]$ \\
$11a_{226}$ & $20/71$ & $4$ & $[4,2,-5,-2]$ & $11a_{229}$ & $16/71$ & $4$ & $[4,-2,3,-2]$ \\
$11a_{230}$ & $8/51$ & $3$ & $[6,-3,-3]$ & $11a_{234}$ & $5/37$ & $3$ & $[7,-2,2]$ \\
$11a_{235}$ & $22/71$ & $4$ & $[3,-4,2,-2]$ & $11a_{236}$ & $29/99$ & $5$ & $[3,-2,2,-2,2]$ \\
$11a_{238}$ & $12/65$ & $4$ & $[5,-2,2,-2]$ & $11a_{242}$ & $9/47$ & $3$ & $[5,-4,2]$ \\
$11a_{243}$ & $20/69$ & $4$ & $[3,-2,4,-2]$ & $11a_{246}$ & $13/41$ & $3$ & $[3,-6,2]$ \\
$11a_{247}$ & $2/19$ & $2$ & $[9,-2]$ & $11a_{306}$ & $29/105$ & $4$ & $[4,3,3,4]$ \\
$11a_{307}$ & $18/83$ & $4$ & $[5,2,-2,-4]$ & $11a_{308}$ & $15/71$ & $3$ & $[5,4,4]$ \\
$11a_{309}$ & $25/93$ & $4$ & $[4,3,-2,-4]$ & $11a_{310}$ & $14/61$ & $3$ & $[4,-3,-5]$ \\
\end{tabular}
\end{center}
\begin{center}
\begin{tabular}{llllllll}
knot & $p/q$ & $\gamma$ & continued fraction & knot & $p/q$ & $\gamma$ & continued fraction \\
\hline
$11a_{311}$ & $18/79$ & $4$ & $[4,-3,-2,3]$ & $11a_{333}$ & $14/65$ & $3$ & $[5,3,5]$ \\
$11a_{334}$ & $9/49$ & $3$ & $[5,-2,4]$ & $11a_{335}$ & $17/75$ & $4$ & $[4,-2,2,-3]$ \\
$11a_{336}$ & $11/59$ & $3$ & $[5,-3,-4]$ & $11a_{337}$ & $26/89$ & $4$ & $[3,-2,3,4]$ \\
$11a_{339}$ & $13/55$ & $3$ & $[4,-4,3]$ & $11a_{341}$ & $19/61$ & $3$ & $[3,-5,-4]$ \\
$11a_{342}$ & $4/29$ & $2$ & $[7,-4]$ & $11a_{343}$ & $4/31$ & $3$ & $[8,4]^*$ \\
$11a_{355}$ & $7/45$ & $3$ & $[6,-2,3]$ & $11a_{356}$ & $24/79$ & $4$ & $[3,-3,2,-3]$ \\
$11a_{357}$ & $27/91$ & $4$ & $[3,-3,-3,3]$ & $11a_{358}$ & $5/31$ & $2$ & $[6,-5]$ \\
$11a_{359}$ & $10/53$ & $3$ & $[5,-3,3]$ & $11a_{360}$ & $10/57$ & $3$ & $[6,3,-3]$ \\
$11a_{363}$ & $6/35$ & $3$ & $[6,6]^*$ & $11a_{364}$ & $3/25$ & $2$ & $[8,-3]$ \\
$11a_{365}$ & $16/51$ & $3$ & $[3,-5,3]$ & $11a_{367}$ & $1/11$ & $1$ & $[11]$ \\
$12a_{38}$ & $33/71$ & $4$ & $[2,-6,2,3]$ & $12a_{169}$ & $23/49$ & $3$ & $[2,-8,-3]$ \\
$12a_{197}$ & $32/69$ & $4$  & $[2,-6,3,2]$  & $12a_{204}$ & $76/173$ & $5$ & $[2,-4,-3,-3,-3]$ \\
$12a_{206}$ & $47/105$ & $4$  & $[2,-4,4,3]$  & $12a_{221}$ & $66/169$ & $5$ & $[3,2,-4,-3,-3]$ \\
$12a_{226}$ & $75/181$ & $6$  & $[2,-2,2,-2,2,3]$  & $12a_{239}$ & $40/87$ & $4$ & $[2,-6,-3,2]$ \\
$12a_{241}$ & $57/127$ & $5$  & $[2,-4,3,2,-2]$  & $12a_{243}$ & $60/133$ & $5$ & $[2,-4,2,3,3]$ \\
$12a_{247}$ & $71/163$ & $5$  & $[2,-3,3,3,3]$  & $12a_{251}$ & $59/159$ & $5$ & $[3,3,-3,2,3]$ \\
$12a_{254}$ & $23/97$ & $4$  & $[4,-4,2,3]$  & $12a_{255}$ & $28/107$ & $4$ & $[4,5,-2,-3]$ \\
$12a_{257}$ & $80/191$ & $6$  & $[2,-2,2,3,-2,-3]$  & $12a_{259}$ & $52/115$ & $4$ & $[2,-5,-4,-3]$ \\
$12a_{300}$ & $68/155$ & $5$  & $[2,-3,2,4,3]$  & $12a_{302}$ & $61/147$ & $5$ & $[2,-2,2,-4,-3]$ \\
$12a_{303}$ & $64/153$ & $5$  & $[2,-2,2,5,3]$  & $12a_{306}$ & $64/147$ & $5$ & $[2,-3,3,3,-2]$ \\
$12a_{307}$ & $69/157$ & $5$  & $[2,-4,-3,-3,2]$  & $12a_{330}$ & $43/95$ & $4$ & $[2,-5,-4,2]$ \\
$12a_{378}$ & $45/127$ & $4$  & $[3,6,3,3]$  & $12a_{379}$ & $17/71$ & $3$ & $[4,-6,-3]$ \\
$12a_{380}$ & $20/77$ & $3$  & $[4,7,3]$  & $12a_{384}$ & $62/151$ & $5$ & $[2,-2,3,-3,-3]$ \\
$12a_{385}$ & $66/161$ & $5$  & $[2,-2,4,3,3]$  & $12a_{406}$ & $74/179$ & $6$ & $[2,-2,3,2,-2,2]$ \\
$12a_{425}$ & $37/81$ & $4$  & $[2,-5,3,-2]$  & $12a_{437}$ & $65/149$ & $5$ & $[2,-3,2,-3,-3]$ \\
$12a_{447}$ & $43/121$ & $4$  & $[3,5,-3,-3]$  & $12a_{454}$ & $27/103$ & $4$ & $[4,5,-2,2]$ \\
$12a_{471}$ & $38/85$ & $4$  & $[2,-4,5,2]$  & $12a_{477}$ & $70/169$ & $6$ & $[2,-2,2,-2,3,2]$ \\
$12a_{482}$ & $22/93$ & $4$  & $[4,-4,3,2]$  & $12a_{497}$ & $81/209$ & $6$ & $[3,2,-3,-2,2,3]$ \\
$12a_{498}$ & $76/207$ & $5$  & $[3,4,3,3,3]$  & $12a_{499}$ & $89/233$ & $6$ & $[3,3,3,3,2,-2]$ \\
$12a_{500}$ & $60/167$ & $5$  & $[3,4,-2,-3,-3]$  & $12a_{501}$ & $55/199$ & $5$ & $[4,3,3,3,3]$ \\
$12a_{502}$ & $37/91$ & $4$  & $[2,-2,6,3]$  & $12a_{506}$ & $68/185$ & $5$ & $[3,3,-2,-4,-3]$ \\
$12a_{508}$ & $56/129$ & $5$  & $[2,-3,3,-2,2]$  & $12a_{510}$ & $81/193$ & $6$ & $[2,-2,2,3,3,-2]$ \\
$12a_{511}$ & $51/125$ & $5$  & $[2,-2,5,2,-2]$  & $12a_{512}$ & $64/151$ & $5$ & $[2,-3,-4,2,3]$ \\
$12a_{514}$ & $79/187$ & $5$  & $[2,-3,-4,-3,-3]$  & $12a_{517}$ & $52/145$ & $4$ & $[3,5,4,3]$ \\
$12a_{518}$ & $34/157$ & $5$  & $[5,2,-2,-3,-3]$  & $12a_{519}$ & $25/111$ & $4$ & $[4,-2,4,3]$ \\
$12a_{520}$ & $36/133$ & $4$  & $[4,3,-4,-3]$  & $12a_{521}$ & $48/113$ & $4$ & $[2,-3,-6,-3]$ \\
$12a_{522}$ & $73/173$ & $5$  & $[2,-3,-3,3,3]$  & $12a_{528}$ & $67/183$ & $5$ & $[3,4,4,2,-2]$ \\
$12a_{532}$ & $33/125$ & $4$  & $[4,5,3,-2]$  & $12a_{533}$ & $31/137$ & $5$ & $[4,-2,3,2,-2]$ \\
$12a_{534}$ & $44/163$ & $5$  & $[4,3,-3,-2,2]$  & $12a_{535}$ & $47/175$ & $5$ & $[4,3,-2,-3,-3]$ \\
$12a_{536}$ & $29/137$ & $4$  & $[5,4,3,3]$  & $12a_{537}$ & $50/179$ & $5$ & $[4,2,-3,-3,-3]$ \\
$12a_{538}$ & $13/83$ & $4$  & $[6,-2,2,3]$  & $12a_{539}$ & $44/145$ & $5$ & $[3,-3,3,2,-2]$ \\
$12a_{540}$ & $49/165$ & $5$  & $[3,-3,-3,2,3]$  & $12a_{541}$ & $41/153$ & $4$ & $[4,4,4,3]$ \\
$12a_{545}$ & $63/143$ & $5$  & $[2,-4,-3,2,-2]$  & $12a_{549}$ & $26/111$ & $4$ & $[4,-4,-3,2]$ \\
$12a_{550}$ & $34/149$ & $5$  & $[4,-2,2,3,3]$  & $12a_{551}$ & $18/103$ & $4$ & $[6,3,-2,-3]$ \\
$12a_{552}$ & $30/131$ & $4$  & $[4,-3,-4,-3]$  & $12a_{579}$ & $49/177$ & $5$ & $[4,2,-2,-4,-3]$ \\
$12a_{580}$ & $11/69$ & $3$  & $[6,-4,-3]$  & $12a_{581}$ & $36/119$ & $4$ & $[3,-3,4,3]$ \\
$12a_{582}$ & $39/131$ & $4$  & $[3,-3,-5,-3]$  & $12a_{583}$ & $45/161$ & $5$ & $[4,2,-3,-3,2]$ \\
$12a_{584}$ & $31/143$ & $5$  & $[5,2,-2,-3,2]$  & $12a_{585}$ & $50/181$ & $5$ & $[4,3,3,3,-2]$ \\
$12a_{595}$ & $30/139$ & $4$  & $[5,3,4,3]$  & $12a_{596}$ & $14/81$ & $3$ & $[6,5,3]$ \\
$12a_{597}$ & $26/123$ & $4$  & $[5,4,3,-2]$  & $12a_{600}$ & $25/109$ & $4$ & $[4,-3,-4,2]$ \\
$12a_{601}$ & $34/127$ & $4$  & $[4,4,5,2]$  & $12a_{643}$ & $23/99$ & $4$ & $[4,-3,3,-2]$ \\
\end{tabular}
\end{center}
\begin{center}
\begin{tabular}{llllllll}
knot & $p/q$ & $\gamma$ & continued fraction & knot & $p/q$ & $\gamma$ & continued fraction \\
\hline
$12a_{644}$ & $30/113$ & $4$  & $[4,4,-3,2]$  & $12a_{649}$ & $27/127$ & $4$ & $[5,3,-3,-3]$ \\
$12a_{650}$ & $46/165$ & $5$  & $[4,2,-2,3,3]$  & $12a_{651}$ & $17/97$ & $4$ & $[6,3,-2,2]$ \\
$12a_{652}$ & $46/155$ & $5$  & $[3,-3,-3,2,-2]$  & $12a_{682}$ & $29/107$ & $4$ & $[4,3,-4,2]$ \\
$12a_{684}$ & $41/135$ & $5$  & $[3,-3,2,-2,2]$  & $12a_{690}$ & $20/89$ & $4$ & $[4,-2,5,2]$ \\
$12a_{691}$ & $12/77$ & $4$  & $[6,-2,3,2]$  & $12a_{713}$ & $39/139$ & $5$ & $[4,2,-3,2,-2]$ \\
$12a_{714}$ & $19/107$ & $4$  & $[6,3,3,-2]$  & $12a_{715}$ & $50/169$ & $5$ & $[3,-3,-3,-3,2]$ \\
$12a_{716}$ & $5/43$ & $3$  & $[9,2,-2]$  & $12a_{717}$ & $28/89$ & $4$ & $[3,-6,-2,2]$ \\
$12a_{718}$ & $41/141$ & $5$  & $[3,-2,4,2,-2]$  & $12a_{720}$ & $21/113$ & $4$ & $[5,-3,-3,-3]$ \\
$12a_{721}$ & $50/171$ & $5$  & $[3,-2,3,3,3]$  & $12a_{722}$ & $3/29$ & $2$ & $[10,3]$ \\
$12a_{723}$ & $20/63$ & $3$  & $[3,-7,-3]$  & $12a_{724}$ & $31/107$ & $4$ & $[3,-2,5,3]$ \\
$12a_{726}$ & $19/103$ & $4$  & $[5,-2,3,3]$  & $12a_{727}$ & $46/157$ & $5$ & $[3,-2,2,-3,-3]$ \\
$12a_{728}$ & $29/133$ & $5$  & $[5,2,-2,2,-2]$  & $12a_{729}$ & $46/167$ & $5$ & $[4,3,3,-2,2]$ \\
$12a_{731}$ & $22/105$ & $4$  & $[5,4,-2,2]$  & $12a_{732}$ & $18/95$ & $4$ & $[5,-3,2,3]$ \\
$12a_{733}$ & $14/73$ & $3$  & $[5,-5,-3]$  & $12a_{736}$ & $43/141$ & $5$ & $[3,-3,2,3,-2]$ \\
$12a_{738}$ & $37/119$ & $4$  & $[3,-5,-3,-3]$  & $12a_{740}$ & $35/113$ & $4$ & $[3,-4,3,3]$ \\
$12a_{743}$ & $12/79$ & $4$  & $[7,2,-2,2]$  & $12a_{744}$ & $8/61$ & $3$ & $[8,3,3]$ \\
$12a_{745}$ & $8/59$ & $3$  & $[7,-3,-3]$  & $12a_{758}$ & $31/113$ & $4$ & $[4,3,5,-2]$ \\
$12a_{759}$ & $9/61$ & $3$  & $[7,4,-2]$  & $12a_{760}$ & $34/111$ & $4$ & $[3,-4,-4,2]$ \\
$12a_{761}$ & $41/139$ & $5$  & $[3,-2,2,4,-2]$  & $12a_{762}$ & $7/51$ & $3$ & $[7,-3,2]$ \\
$12a_{763}$ & $30/97$ & $4$  & $[3,-4,3,-2]$  & $12a_{764}$ & $39/133$ & $5$ & $[3,-2,2,-3,2]$ \\
$12a_{773}$ & $20/91$ & $4$  & $[4,-2,-5,2]$  & $12a_{774}$ & $16/89$ & $4$ & $[5,-2,-4,2]$ \\
$12a_{775}$ & $16/87$ & $4$  & $[5,-2,3,-2]$  & $12a_{791}$ & $13/63$ & $3$ & $[5,6,-2]$ \\
$12a_{792}$ & $24/85$ & $4$  & $[4,2,-5,2]$  & $12a_{796}$ & $11/57$ & $3$ & $[5,-5,2]$ \\
$12a_{797}$ & $24/83$ & $4$  & $[3,-2,5,-2]$  & $12a_{802}$ & $15/47$ & $3$ & $[3,-7,2]$ \\
$12a_{803}$ & $2/21$ & $2$  & $[11,2]$  & $12a_{1023}$ & $29/127$ & $4$ & $[4,-3,-3,-4]$ \\
$12a_{1024}$ & $40/149$ & $4$  & $[4,4,3,4]$  & $12a_{1029}$ & $19/81$ & $3$ & $[4,-4,-5]$ \\
$12a_{1030}$ & $19/91$ & $3$  & $[5,5,4]$  & $12a_{1033}$ & $25/107$ & $4$ & $[4,-3,2,4]$ \\
$12a_{1034}$ & $32/121$ & $4$  & $[4,5,2,-3]$  & $12a_{1039}$ & $37/137$ & $4$ & $[4,3,-3,-4]$ \\
$12a_{1040}$ & $26/115$ & $4$  & $[4,-2,3,4]$  & $12a_{1125}$ & $23/101$ & $4$ & $[4,-2,2,5]$ \\
$12a_{1126}$ & $26/119$ & $4$  & $[5,2,-3,-4]$  & $12a_{1127}$ & $22/97$ & $4$ & $[5,2,3,-4]$ \\
$12a_{1128}$ & $9/59$ & $3$  & $[7,2,-4]$  & $12a_{1129}$ & $23/105$ & $4$ & $[4,-2,-4,3]$ \\
$12a_{1130}$ & $27/125$ & $4$  & $[5,3,3,-3]$  & $12a_{1131}$ & $11/73$ & $3$ & $[7,3,4]$ \\
$12a_{1132}$ & $40/131$ & $4$  & $[3,-4,-3,-4]$  & $12a_{1133}$ & $47/159$ & $5$ & $[3,-2,2,3,4]$ \\
$12a_{1134}$ & $7/53$ & $3$  & $[8,2,-3]$  & $12a_{1135}$ & $32/103$ & $4$ & $[3,-4,2,4]$ \\
$12a_{1136}$ & $43/147$ & $5$  & $[3,-2,3,2,-3]$  & $12a_{1138}$ & $14/79$ & $3$ & $[6,3,5]$ \\
$12a_{1139}$ & $18/101$ & $4$  & $[6,3,2,-3]$  & $12a_{1140}$ & $18/97$ & $4$ & $[5,-2,2,4]$ \\
$12a_{1145}$ & $15/79$ & $3$  & $[5,-4,-4]$  & $12a_{1146}$ & $34/117$ & $4$ & $[3,-2,4,4]$ \\
$12a_{1148}$ & $23/73$ & $3$  & $[3,-6,-4]$  & $12a_{1149}$ & $4/35$ & $2$ & $[9,4]$ \\
$12a_{1157}$ & $5/39$ & $2$  & $[8,5]$  & $12a_{1158}$ & $16/77$ & $3$ & $[5,5,-3]$ \\
$12a_{1159}$ & $24/113$ & $4$  & $[5,3,-2,3]$  & $12a_{1161}$ & $14/75$ & $3$ & $[5,-3,-5]$ \\
$12a_{1162}$ & $13/69$ & $3$  & $[5,-3,4]$  & $12a_{1163}$ & $24/103$ & $4$ & $[4,-3,2,-3]$ \\
$12a_{1165}$ & $16/67$ & $3$  & $[4,-5,3]$  & $12a_{1166}$ & $4/33$ & $3$ & $[8,-4]^*$ \\
$12a_{1273}$ & $11/61$ & $3$  & $[6,2,-5]$  & $12a_{1274}$ & $17/95$ & $4$ & $[6,2,-2,3]$ \\
$12a_{1275}$ & $44/149$ & $5$  & $[3,-3,-2,2,-3]$  & $12a_{1276}$ & $13/75$ & $3$ & $[6,4,-3]$ \\
$12a_{1277}$ & $37/121$ & $4$  & $[3,-4,-3,3]$  & $12a_{1278}$ & $6/41$ & $2$ & $[7,6]$ \\
$12a_{1279}$ & $10/67$ & $3$  & $[7,3,-3]$  & $12a_{1281}$ & $33/109$ & $4$ & $[3,-3,3,-3]$ \\
$12a_{1282}$ & $10/63$ & $3$  & $[6,-3,3]$  & $12a_{1287}$ & $6/37$ & $3$ & $[6,-6]^*$ \\
\end{tabular}
\end{center}

\bibliographystyle{amsplain}

\end{document}